\documentclass[11pt,bezier]{article}
\usepackage{amsmath,amssymb,amsfonts,euscript,pifont,mathrsfs,latexsym}

\newtheorem{prethm}{{\bf Theorem}}

\newenvironment{thm}{\begin{prethm}\sl{\hspace{-0.5
               em}{\bf.}}}{\end{prethm}}

\newtheorem{prepro}[prethm]{{\bf Proposition}}

\newtheorem{prelem}[prethm]{{\bf Lemma}}

\newenvironment{lem}{\begin{prelem}\sl{\hspace{-0.5
               em}{\bf.}}}{\end{prelem}}

\newtheorem{precor}[prethm]{{\bf Corollary}}

\newtheorem{preconj}[prethm]{{\bf Conjecture}}

\newtheorem{preremark}[prethm]{{\bf Remark}}

\newtheorem{preexample}[prethm]{{\bf Example}}

\newtheorem{prethmm}{{\bf Theorem}}

\newtheorem{preproof}{{\bf{Proof.}}}

\newenvironment{proof}[1]{\begin{preproof}{\rm
               #1}\hfill{$\Box$}}{\end{preproof}}

\DeclareMathAlphabet{\mathpzc}{OT1}{pzc}{m}{it}

\title{\bf\LARGE  Integral trees of  odd  diameters}

\author{\large E. Ghorbani \quad \quad  A. Mohammadian \quad \quad B. Tayfeh-Rezaie\thanks{Corresponding author.}\\[.3cm]
{\sl School of Mathematics, Institute for Research in Fundamental
Sciences (IPM),}\\{\sl P.O. Box
19395-5746, Tehran, Iran }
\\[.3cm]{
$\mathsf{e\_ghorbani@ipm.ir}$ \quad\quad  $\mathsf{ali\_m@ipm.ir}$ \quad\quad  $\mathsf{tayfeh}$-$\mathsf{r@ipm.ir}$}}

\date{}

\begin{document}
\maketitle

\vspace{5mm}

\begin{abstract}
\noindent
A graph is called integral if all eigenvalues of its adjacency matrix consist entirely of integers. Recently,  Csikv\'{a}ri proved the existence of integral  trees
of any   even  diameter.  In the odd case, integral  trees have been constructed with diameter at most $7$. In this paper, we show  that for every odd integer  $n>1$,  there are infinitely many integral trees of diameter  $n$.

\vspace{5mm}
\end{abstract}

\vspace{5mm}

\section{Introduction}

Let $G$ be a graph with the vertex set $\{v_1, \ldots , v_n\}$. The {\sl adjacency matrix} of $G$ is an $n \times  n$
matrix $A(G)$ whose $(i, j)$-entry is $1$ if $v_i$ is adjacent to $v_j$ and  $0$ otherwise.  The   {\sl characteristic polynomial} of  $G$,   denoted by $\varphi(G; x)$, is   the characteristic polynomial of $A(G)$. We will drop the  indeterminate  $x$  for simplicity of notation whenever there is no danger  of  confusion. The zeros of $\varphi(G)$
are called the {\sl eigenvalues} of $G$. Note that  $A(G)$ is a real symmetric matrix so that all eigenvalues of $G$ are reals. The graph $G$ is  said to be
{\sl integral} if all  eigenvalues of $G$  are integers.

The notion of integral graphs was f
irst introduced in   \cite{har}. The  general characterization   problem of
 integral graphs seems to be intractable. Therefore,  it is natural to deal with the problem within specific classes of graphs such as trees, cubic graphs and so on.  Here, we are concerned  with   finding integral trees.  These objects are  extremely rare and very  difficult to find.  For instance, among trees up to $50$ vertices  there are only $28$ integral ones   \cite{b} and   out of a total number  of  $2,262,366,343,746$ trees on $35$ vertices only one tree is integral.
For a long time,   it has been an open question  whether there exist integral trees of any   diameter   \cite{watsc}. Many   attempts by various authors
 led to the constructions of integral trees of diameters $2$--$8$ and  $10$, see   \cite{{hic}, {hic1}, {hic2}, {liu}, {wan}, {wan1}, {wat}}.
 Recently,  Csikv\'{a}ri  found   a nice   construction of  integral trees for  any  even diameter   \cite{csi}.
In the odd case, all integral trees of diameter  $3$ have been characterized by a result from number theory given in \cite{gra}. Infinitely
many integral trees of diameter $5$ were first constructed in \cite{liu}. The existence  of integral trees of diameter $7$ was established  in \cite{hic1} where  the authors found  four such  trees.
In this paper, we show  that for every odd integer  $n>1$,  there are infinitely many integral trees of diameter $n$.

\section{Csikv\'{a}ri's  trees}

In this section,  we revisit the trees constructed in   \cite{csi} and compute their eigenvalues  using   a simple argument.  We will make use of them  to construct new integral trees of odd diameters.

A {\sl rooted tree} is a tree with  a specified  vertex called the  {\sl root}. For a  rooted tree $T$, we denote by $T'$
  the forest resulting from removing the root of  $T$.
 Let    $n$ be a positive integer and  $T_1, T_2$ be two  rooted trees  with  disjoint vertex sets.  Then  $T_1\thicksim nT_2$ is the rooted tree obtained from
 $T_1$ and $n$ copies of $T_2$  by joining the root of $T_1$ to the roots of the copies of $T_2$. The root of $T_1$ is considered
to  be  the root  of the  resulting tree.
For    positive integers $r_1<r_2< \cdots<r_n$,  we define  the rooted tree  $C(r_1, \ldots  , r_n)$  constructed in   \cite{csi}  recursively as
$$
  C(r_1, \ldots  , r_n)=C(r_1, \ldots  , r_{n-2})\thicksim (r_n-r_{n-1})C(r_1, \ldots  , r_{n-1}),
$$
for  $n\geqslant2$,  with initial trees $C(\,)$ and $C(r_1)$ being  the   one-vertex  tree and the  star tree on $r_1+1$ vertices, respectively.

The following lemma is proved in   \cite[p.\,59]{cve} for $n=1$. The general case is straightforward by induction on $n$.

\begin{lem}\label{main1}
Let  $T_1, T_2$ be two  rooted  trees and let $T = T_1\thicksim nT_2$.  Then
$$\varphi(T )=\varphi(T_2 )^{n-1}\left(\varphi(T_1 )\varphi(T_2 )-n\varphi(T_1')\varphi(T_2')\right).$$
\end{lem}

The following lemma can be used to determine  the eigenvalues of Csikv\'{a}ri's  trees and  their multiplicities.

\begin{lem}\label{csi}
Let   $n\geqslant2$ and  $r_1, \ldots , r_n$ be positive integers. Then
$$\varphi\big(C(r_1, \ldots, r_n)\big)=\varphi^{r_n-r_{n-1}}\big(C(r_1, \ldots, r_{n-1})\big)\varphi\big(C(r_1, \ldots, r_{n-2})\big)\frac{x^2-r_n}{x^2-r_{n-1}}.$$
\end{lem}

\begin{proof}{For convenience,  we set  $P_k=\varphi(C(r_1, \ldots, r_k))$, $Q_k=\varphi(C'(r_1, \ldots, r_k))$  and $d_k=r_k-r_{k-1}$ for $k\geqslant1$, with the
convention $r_0=0$.
Since
$C'(r_1, \ldots  , r_k)=C'(r_1, \ldots  , r_{k-2})\cup d_kC(r_1, \ldots  , r_{k-1})$
 for  $k\geqslant2$, we have $Q_k=P_{k-1}^{d_{k}}Q_{k-2}$. By Lemma   \ref{main1},
\begin{align*}
P_k&=P^{d_k-1}_{k-1}\big(P_{k-1}P_{k-2}-(r_k-r_{k-1})Q_{k-1}Q_{k-2}\big)\nonumber\\
&= P^{d_k-1}_{k-1}\left(P^{d_{k-1}}_{k-2}\big(P_{k-2}P_{k-3}-(r_{k-1}-r_{k-2})Q_{k-2}Q_{k-3}\big)-(r_k-r_{k-1})P_{k-2}^{d_{k-1}}Q_{k-2}Q_{k-3}\right)
\nonumber\\
&= P^{d_k-1}_{k-1}P^{d_{k-1}}_{k-2}\big(P_{k-2}P_{k-3}-(r_k-r_{k-2})Q_{k-2}Q_{k-3}\big)\nonumber\\
&\,\,\,\vdots\nonumber\\
&=P^{d_k-1}_{k-1}P^{d_{k-1}}_{k-2}\cdots P_2^{d_3}\big(P_2P_1-(r_k-r_{2})Q_2Q_1\big)\nonumber\\
&=P^{d_k-1}_{k-1}P^{d_{k-1}}_{k-2}\cdots P_1^{d_2}\big(P_1x-(r_k-r_{1})Q_1\big)\nonumber\\
&=P^{d_k-1}_{k-1}P^{d_{k-1}}_{k-2}\cdots P_1^{d_2}x^{d_1}(x^2-r_k).
\end{align*}
Note that    $P_1=x^{d_1-1}(x^2-r_1)$ and so $P_k=P^{d_k-1}_{k-1}P^{d_{k-1}}_{k-2}\cdots P_1^{d_2}x^{d_1}(x^2-r_k)$  holds for  $k\geqslant1$.  To complete the proof, apply this equality    for  $k=n-1, n$ and  then    compute $P_n/P_{n-1}$.
}\end{proof}

 It is not hard to see that  diameter of  $C(r_1, \ldots, r_n)$ is $2n$ provided that $r_n-r_{n-1}>1$. The following theorem readily  follows from Lemma   \ref{csi} which     establishes  the existence of  infinitely many  integral trees of any even   diameter.

\begin{thm} \label{csith}
{\rm\cite{csi}}
The set of distinct  eigenvalues of the tree $C(r_1, \ldots, r_n)$ is $\left\{0, \pm\sqrt{r_1},  \ldots, \pm\sqrt{r_n}\right\}$.
\end{thm}

Let us  introduce an alternative   representation of
$\varphi(C(r_1, \ldots, r_n))$ and $\varphi(C'(r_1, \ldots, r_n))$
 to be used  in  sequel.
For  $C=C(r_1, \ldots, r_n)$,  we let
$$f(C)=\prod_{i=1}^{\left \lceil\frac{n}{2}\right\rceil}\frac{\varphi^{d_{n-2i+2}}\big(C(r_1, \ldots, r_{n-2i+1})\big)}{x^2-r_{n-2i+1}},$$
where $d_i=r_{i}-r_{i-1}$ with the convention $r_0=0$. By Lemma   \ref{csi}, $f(C)$ is a polynomial and  clearly   we  have
\begin{equation}\label{phif}
\varphi(C)=xf(C)\prod_{i=1}^{\left\lceil \frac{n}{2}\right\rceil}(x^2-r_{n-2i+2})
\end{equation} and
\begin{equation}\label{qf}
\varphi(C')=f(C)\prod_{i=1}^{ \left\lceil\frac{n}{2}\right\rceil}(x^2-r_{n-2i+1}).
\end{equation}
Note that by Lemma  \ref{csi}, if $r_n-r_{n-1}>1$, then the positive zeroes of $f(C)$ read as
$\sqrt{r_1},  \ldots, \sqrt{r_{n-1}}.$

\section{A  class of  trees}\label{3}

In this section,  we introduce  a class of trees which will be  used  to obtain integral trees of odd diameters.
For   positive integers  $n, r, r_0, r_1,   \ldots, r_n$   such that  $n\geqslant2$ and  $\max\{r_0, r_1\}<r_2<\cdots<r_n$,
let $U=C(r_1, \ldots, r_n)$, $V=C(r_0,r_2, \ldots, r_{n-1})$, $W=C(r_2,\ldots, r_n)$,  and define
$$T(r, r_0,r_1, \ldots,r_n)=U\thicksim (V\thicksim rW).$$
Note that for $n=2$, we let $V=C(r_0)$.
It is easily checked that  the  maximum   distance between  a vertex   of      $C(k_1, \ldots,k_n)$  and  its    root is $n$. So
 $T=T(r, r_0, r_1, \ldots, r_n)$ is a tree  of diameter $2n+1$.
We   proceed to   determine   $\varphi(T)$.
Applying  Lemma \ref{main1}, we find that
\begin{align*}
\varphi(T)=&\,\varphi(U)\varphi^{r-1}(W)\left(\varphi(V)\varphi(W)-r\varphi(V')\varphi(W')\right)-\varphi(U')\varphi(V')\varphi^r(W)\nonumber\\
=&\varphi^{r-1}(W)\left(\varphi(U)\varphi(V)\varphi(W)-r\varphi(U)\varphi(V')\varphi(W')-\varphi(U')\varphi(V')\varphi(W)\right).
\end{align*}

First assume that $n=2m+1$ is odd.
By  (\ref{phif}) and (\ref{qf}), we have
\begin{align}
\varphi(U)=&\,xf(U)(x^2-r_1)(x^2-r_n)\prod_{i=2}^{m}(x^2-r_{2i-1}),\label{odd1}\\
\varphi(V)=&\,xf(V)\prod_{i=1}^{m}(x^2-r_{2i}),\\
\varphi(W)=&\,xf(W)(x^2-r_n)\prod_{i=2}^{m}(x^2-r_{2i-1}),
\end{align}
and
\begin{align}
\varphi(U')=&\,x^2f(U)\prod_{i=1}^{m}(x^2-r_{2i}),  \\
 \varphi(V')=&\,f(V)(x^2-r_0)\prod_{i=2}^{m}(x^2-r_{2i-1}),\\
 \varphi(W')=&\,f(W)\prod_{i=1}^{m}(x^2-r_{2i}).\label{odd2}
\end{align}
Hence, by  (\ref{odd1})--(\ref{odd2}),
\begin{equation*}\label{char}
\varphi(T)=x(x^2-r_n)\varphi^{r-1}(W)f(U)f(V)f(W)\prod_{i=2}^{m}(x^2-r_{2i})\prod_{i=2}^{m}(x^2-r_{2i-1})^2\psi_o(x),
\end{equation*}
where
\begin{equation*}
\psi_o(T)=x^2(x^2-r_1)(x^2-r_n)-r(x^2-r_0)(x^2-r_1)-x^2(x^2-r_{0}).
\end{equation*}

Next  suppose  that $n=2m$ is even.
By  (\ref{phif}) and (\ref{qf}), we have
\begin{align}
\varphi(U)=&\,xf(U)(x^2-r_n)\prod_{i=1}^{m-1}(x^2-r_{2i}),\label{even1}\\
\varphi(V)=&\,xf(V)(x^2-r_0)\prod_{i=2}^{m}(x^2-r_{2i-1}),\\
\varphi(W)=&\,xf(W)(x^2-r_n)\prod_{i=1}^{m-1}(x^2-r_{2i}),
\end{align}
and
\begin{align}
\varphi(U')=&\,f(U)(x^2-r_{1})\prod_{i=2}^{m}(x^2-r_{2i-1}),  \\
 \varphi(V')=&\,x^2f(V)\prod_{i=1}^{m-1}(x^2-r_{2i}),\\
 \varphi(W')=&\,x^2f(W)\prod_{i=2}^{m}(x^2-r_{2i-1}).\label{even2}
\end{align}
Thus, by  (\ref{even1})--(\ref{even2}),
\begin{equation*}\label{char}
\varphi(T)=x^3(x^2-r_n)\varphi^{r-1}(W)f(U)f(V)f(W)\prod_{i=2}^{m}(x^2-r_{2i-1})\prod_{i=1}^{m-1}(x^2-r_{2i})^2\psi_e(x),
\end{equation*}
where
\begin{equation*}
\psi_e(T)=(x^2-r_0)(x^2-r_n)-rx^2-(x^2-r_1).
\end{equation*}

In summary, using the above notation, we have the following theorem.
\begin{thm}\label{intT}
Let $n$ be odd (respectively, even). Then  $T$ is an integral tree of diameter $2n+1$ if and only if $r_0,r_1,\ldots,r_n$ are
perfect squares and all the zeros of  $\psi_o(T)$ (respectively, $\psi_e(T)$) are integers.
\end{thm}

\section{Integral trees of diameter $4k+1$}

Let $n$ be  even.
It is not difficult to choose the   parameters of $T$ in such a way  that  the zeros of  $\psi_e(T)$ are all integers.
For instance, let $r_0=1$, $r_1=4k^2$, $r_n=(k^2-1)^2$ and  $r=4k^2-1$.
Then
\begin{align*}
\psi_e(T)=&\,(x^2-1)\big(x^2-(k^2-1)^2\big)-(4k^2-1)x^2-(x^2-4k^2)\\
=&\,(x^2-1)\big(x^2-(k^2+1)^2\big).
\end{align*}
Clearly,  if we   choose $k$ large enough, then we are able to   take distinct   prefect squares $r_2, \ldots, r_{n-1}$ in the interval $(4k^2, (k^2-1)^2)$.   Hence, we have proved  the  following theorem.

\begin{thm}
For every even positive integer $n$,  there are infinitely many integral trees of diameter  $2n+1$.
\end{thm}

\section{Integral trees of diameter $4k+3$}

Let $n$ be  odd.
Our   objective   is to  choose  the parameters of $T$  in such a way that  all  the zeros of   $\psi_o(T)$ are  integers. This can be done in many ways. For instance,  if we set $r_0=r_1=a^2$ and $r=r_n=4(a-1)^2$ for some integer $a$ with  $|a|\geqslant3$, then
$$\psi_o(T)=\left(x^2-a^2\right)\Big(x^4-\big(8(a-1)^2+1\big)x^2+4a^2(a-1)^2\Big).$$
The zeros of $\psi_o(T)$ are $\pm a$ and $\pm(a-\tfrac{3}{2})\pm\tfrac{1}{2}\sqrt{12a^2-20a+9}$.
So the zeros of $\psi_o(T)$ are    integers if and only if  $12a^2-20a+9$ is a perfect square, say $b^2$. We  have    $(6a-5)^2-3b^2=-2$. From number theory, we know  that the Pell-like equation $x^2-3y^2=-2$ has infinitely many integral solutions with   $x\equiv\pm1\pmod6$. For example, one  may take
$$a=\frac{1}{12}\left(\big(1-\sqrt{3}\big)\big(-2+\sqrt{3}\big)^k+\big(1+\sqrt{3}\big)\big(-2-\sqrt{3}\big)^k+10\right),$$ for arbitrary integer $k\geqslant2$.
Therefore,    we come up with   the  following theorem.
\begin{thm}
For every odd integer  $n\geqslant3$,  there are infinitely many integral trees of diameter  $2n+1$.
\end{thm}

\section*{Acknowledgments}
The research of the first and  second authors was in part supported by a grant from IPM.

\bibliographystyle{siam}

\end{document}